\documentclass[11pt]{article}
\usepackage{amsmath,amssymb,theorem}
\usepackage{graphicx, enumerate}
\usepackage{subfigure}
\usepackage{psfrag}
\usepackage{mathrsfs}
\usepackage{placeins}
\usepackage{wrapfig}
\usepackage{booktabs}
\usepackage{chngcntr}
\usepackage{epstopdf}
\counterwithin{figure}{section}
\numberwithin{figure}{section}
\counterwithin{table}{section}
\usepackage{hyperref}
\hypersetup{
  colorlinks,
  citecolor=black,
  filecolor=black,
  linkcolor=black,
  urlcolor=black
}

\newtheorem{thm}{Theorem}%[section]
\newtheorem{conj}{Conjecture}
\newtheorem{cor}{Corollary}%[section]
\newtheorem{lem}{Lemma}%[section]
\theorembodyfont{\rmfamily}
\newtheorem{defn}{Definition}%[section]
\def\pf{\bigskip\noindent {\bf Proof.}~~}

\def\less{\backslash}

\def\pf{\bigskip\noindent {\bf{Proof.}}~~}

\newcounter{counter}

\def\proofsquare{  \bigskip\hfill\vrule height3pt width6pt depth2pt}
\begin{document}
\title{On Lagrangians of $3$-uniform hypergraphs}
\author{
Hui Lei
\thanks{Center for Combinatorics and LPMC, Nankai University, Tianjin 300071, China,
   ({\tt hlei@mail.nankai.edu.cn}). This author was supported in part by National Natural Science Foundation of China (No. 11771221).}
%  ({\tt leihui0711@163.com}).}
\and
Linyuan Lu
\thanks{University of South Carolina, Columbia, SC 29208,
({\tt lu@math.sc.edu}). This author was supported in part by NSF
grant DMS 1600811.}
\and
Yuejian Peng
\thanks{Hunan University, Changsha, China, ({\tt ypeng1@hnu.edu.cn}). This author was supported in part by National Natural Science Foundation of China (No. 11671124).
}
}
\maketitle

\begin{abstract}
  Frankl and F\"uredi conjectured in 1989 that the maximum Lagrangian of all $r$-uniform
hypergraphs of fixed size $m$ is realized by the minimum hypergraph $C_{r,m}$ under the colexicographic
order. In this paper, we prove a weaker version of the Frankl and F\"{u}redi's conjecture at $r=3$:
there exists an absolute constant $c>0$ such that for any $3$-uniform hypergraph $H$ with $m$ edges,
the Lagrangian of $H$ satisfies $\lambda(H)\leq \lambda(C_{3,m+cm^{2/9}})$.

In particular, this result implies that the Frankl and F\"{u}redi's conjecture holds for $r=3$ and $m\in [{t-1\choose 3}, {t\choose 3}-(t-2)-ct^{\frac{2}{3}}]$. It improves a recent result of Tyomkyn.
\end{abstract}
%{\it{Keywords}}: \\
%{\it {2010 Mathematics Subject Classification}}: 05C55;  05D10; 05C15
\section{Introduction}

For a set $V$ and a positive integer $r$, let $V^{(r)}$ be the family of all $r$-subsets of $V$. An {\it $r$-uniform hypergraph} $G$, (or {\it $r$-graph}, for short), consists of a set $V$ of vertices and a set $E\subseteq V^{(r)}$ of edges. For
an integer $n\in \mathbb{N}$, we denote the set $\{1, 2, 3,\dots, n\}$ by $[n]$. Let
$K_t^{(r)}$ (or $[t]^{(r)}$) denote the {\it complete $r$-graph of order $t$}, that is, the $r$-graph of order $t$ containing all possible edges. Given an $r$-graph $G$, we use $e(G)$ to denote the number of edges of $G$.

%An edge
%$e=\{a_1,a_2,\dots,a_r\}$ will be denoted by $a_1a_2\dots a_r$. For
%an integer $n\in \mathbb{N}$, we denote the set $\{1, 2, 3,\dots, n\}$ by $[n]$. Let
%$K_t^{(r)}$ (or $[t]^{(r)}$) denote the {\it complete $r$-graph of order $t$}, that is, the $r$-graph of order $t$ containing all possible edges.  Given an $r$-graph $G$, we use $|G|$ to denote the number of vertices and $e(G)$ the number of edges of $G$, respectively. The clique number of an $r$-graph $G$, denoted as $w(G)$, is defined as the
%cardinality of the maximum complete $r$-graph in $G$.
\medskip
\begin{defn}
For an $r$-graph
$G$ of order $n$ and a vector $\overset{\rightarrow}{x}=(x_1,\dots,x_n)\in \mathcal{R}^n$, the {\it weight polynomial} of $G$ is
$$w(G,\overset{\rightarrow}{x})=\sum_{e\in E}\prod_{i\in e}x_i.$$
\end{defn}
\begin{defn}
We call $\overset{\rightarrow}{x}=(x_1,\dots,x_n)\in \mathcal{R}^n$  a {\it legal weighting} for $G$ if $x_i\geq0$ for any $i\in[n]$ and $\sum_{i=1}^n x_i=1$.
% \begin{description}
%   \item[(i)] $\forall i\in[n]$~~~  $x_i\geq0$,
%   \item[(ii)]  $\sum_{i=1}^n x_i=1.$
% \end{description}
 \end{defn}

\begin{defn}
The {\it Lagrangian} of $G$ is defined to be $\lambda(G) = \max w(G, \overset{\rightarrow}{x})$, where the maximum is
over all legal weightings for $G$.  We call a legal weighting $\overset{\rightarrow}{x}$ {\it optimal} if $ w(G, \overset{\rightarrow}{x})=\lambda(G)$.
%\begin{description}
%\item[(iii)] $ w(G, \overset{\rightarrow}{x})=\lambda(G)$.
%\end{description}
\end{defn}

Lagrangians for graphs (i.e, $2$-graphs) were introduced by Motzkin and Straus in 1965 \cite{MS1965}. They determined the following simple expression for the Lagrangian of a graph.
\begin{thm}[\cite{MS1965}]\label{graph}
If $G$ is a graph in which a largest clique has order $t$, then
$$\lambda(G)=\lambda(K_t^{(2)})=\frac{1}{2}(1-\frac{1}{t}).$$
\end{thm}
This theorem implies Tur\'an theorem; and Lagrangians are closely related to Tur\'an densities.

Let
\begin{equation}
   \label{eq:mum}
\lambda_r(m)=\max\{ \lambda(H)\colon H \mbox{ is an $r$-graph with $m$ edges}\}.
 \end{equation}
There are rich literatures on determining/estimating the values of $\lambda_r(m)$.

For distinct $A,B\in \mathbb{N}^{(r)}$, we say that $A$ is less than $B$  in the colexicographic ordering if $\max(A\bigtriangleup B)\in B$.
Let $C_{r,m}$ be the subgraph of $\mathbb{N}^{(r)}$ consisting of the first $m$ sets in the colexicographic ordering. If $r=3$, we simply write $C_m$ instead of $C_{3,m}$.

In 1989, Frankl and F\"{u}redi \cite{FF1989} made the following conjecture.
\begin{conj}[\cite{FF1989}]\label{conjecture}
  For any $r\geq 3$ and $m\geq 1$, we have $\lambda_r(m)=\lambda(C_{r,m})$.
 % For any $r\geq 3$ and $m\geq 1$, the maximum Lagrangian of an $r$-graph $H$ with $m$ edges is achieved by $C_{r,m}$;
%i.e.,
%$$\lambda(C_{r,m})=\max\{\lambda(G):G\subseteq \mathbb{N}^{(r)}, e(G)=m\}.$$
\end{conj}

% Lagrangians for graphs (i.e, $2$-graphs) were introduced by Motzkin and Straus in 1965 \cite{MS1965}. They determined the following simple expression for the Lagrangian of a graph.
% \begin{thm}[\cite{MS1965}]\label{graph}
% If $G$ is a graph in which a largest clique has order $t$, then
% $$\lambda(G)=\lambda(K_t^{(2)})=\frac{1}{2}(1-\frac{1}{t}).$$
% \end{thm}

For $r=2$, the validity of Conjecture \ref{conjecture} follows from Theorem \ref{graph}. However, this conjecture is still open even for the first case $r=3$.

Talbot \cite{T2002} has shown that $\lambda(C_{r,m})$ is a constant ($={t-1\choose r}/(t-1)^r$) for
$m \in [{t-1\choose r}, {t\choose r}-{t-2\choose r-2}]$ and jumps for every $m\in [{t\choose r}-{t-2 \choose r-2}, {t\choose r}]$. Most known results are in the interval $[{t-1\choose r}, {t \choose r}-{t-2\choose r-2}]$.  For $r=3$, Talbot \cite{T2002} first proved that Conjecture \ref{conjecture} holds whenever ${{t-1}\choose 3}\leq m\leq {t \choose 3}-{{t-2}\choose 1}-(t-1)={t\choose3}-(2t-3)$ for some $t\in \mathbb{N}$. Tang, Peng, Zhang and Zhao \cite{TPZZ2016}
extended the above range to ${{t-1}\choose 3}\leq m\leq {t\choose 3}-{{t-2}\choose 1}-\frac 1 2(t-1)$. Recently, Tyomkyn \cite{T2017} proved the following.
\begin{thm}[\cite{T2017}]\label{3/4}
  \begin{enumerate}
  \item For $r=3$,  there exists a constant $\delta_3>0$ such that
    for any $m$ satisfying ${{t-1}\choose3}\leq m\leq {t\choose 3}- {t-2\choose 1}-\delta_3t^{3/4}$  we have $$\lambda_3(m)=\frac{{t-1\choose 3}}{(t-1)^3}.$$
  \item For $r\geq 4$,  there exists a constant $\delta_r>0$ such that
    for any $m$ satisfying ${{t-1}\choose r}\leq m\leq {t\choose r}-\delta_rt^{r-2}$ we have
    $$\lambda_r(m)=\frac{{t-1\choose r}}{(t-1)^r}.$$
  \end{enumerate}
\end{thm}

A few good upper bounds on $\lambda(G)$ are known for general $m$. The following result,
which was conjectured (and partially solved for $r=3,4,5$ and any $m$; and for the case $m\geq {4(r-1)(r-2)\choose r}$) by Nikiforov \cite{Nikiforov2018},
was completely proved by the second author.

\begin{thm}[{Lu\cite{maxPspec}}] \label{smooth}
 Let $r \geq 2$ and $H$ be an $r$-uniform hypergraph with $m$ edges. Write $m = {s\choose r}$
   for some real $s\geq r-1$. We have $$\lambda(H)\leq m s^{-r}.$$
   The equality holds if and only if $s$ is an integer and $H$ is the complete  $r$-uniform hypergraph
   $K^r_s$  possibly with some isolated vertices added.
 \end{thm}

 \begin{figure}[hbt]
\centering
\includegraphics[width=300pt, height=200pt]{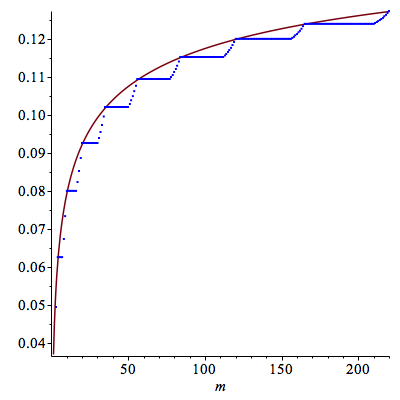}
\caption{The conjectured values of $\lambda_3(m)$ and its smooth upper bound in Theorem \ref{smooth}.}
\label{fig:1}
\end{figure}

The Lagrangians of $3$-graphs have been extensively studied.
In this paper, we focus on $3$-graphs. We would like to prove a better upper bound for $\lambda_3(m)$.
We have the following theorem.
\begin{thm}\label{mainthm}
  There exists a constant $c>0$ such that for any $m>0$  we have
  \begin{equation}
    \label{eq:1}
    \lambda_3(m)\leq \lambda(C_{m+cm^{2/9}}).
  \end{equation}
\end{thm}
Compared with the result of Theorem \ref{smooth} at $r=3$, the upper bound in Theorem \ref{mainthm} is better for most values $m$.
Note that $\lambda(C_m)={t-1\choose 3}/(t-1)^3$ for ${{t-1}\choose3}\leq m\leq {t\choose 3}-{t-2\choose 1}$.
We have the following corollary, which improves Tyomkyn's result for $r=3$ (Theorem \ref{3/4} item 1).
\begin{cor}\label{maincor}
  There exists a constant $c>0$ such that for any ${{t-1}\choose3}\leq m\leq {t\choose3}-(t-2)-ct^{2/3}$,
  we have $$\lambda_3(m)=\frac{{t-1\choose 3}}{(t-1)^3}.$$
%$$\lambda(C_{r,m})=\max\{\lambda(G):G\subseteq \mathbb{N}^{(r)}, e(G)=m\}.$$
\end{cor}

% {\bf Remark:} The proof of Theorem \ref{mainthm} actually gives a better upper bound
% when $m$ is close to ${t\choose 3}$. Let $m={t\choose 3}-l$ where $0<l<(t-2)+ct^{2/3}$.
% Then we actually prove
% $$\lambda_3(m)\leq \lambda(C_{m+cl^{2/3}}).$$

The paper is organized as follows. In section 2, we review notation and facts. Theorem \ref{mainthm} will be proved in section 3.

\section{Notation and Preliminaries}

Although our paper is focusing on $r=3$, we would like to give preliminaries for general $r$ first.
\subsection{General $r$}
Let $r\geq2$ be an integer. Given an $r$-graph $G=(V,E)$ and a set $S\subseteq \mathbb{N}$ with $|S|<r$, the $(r-|S|)$-uniform {\it link hypergraph} of $S$   is defined as $G_S=(V,E_S)$ with $E_S:=\{A\in \mathbb{N}^{(r-|S|)}:A\cup S\in E\}$. We will denote the complement graph of $G_S$ by $G^c_S=(V,E^c_S)$ with $E^c_S:=\{A\in \mathbb{N}^{(r-|S|)}:A\cup S\in V^{(r)}\less E\}$. Define $G_{i\less j}=(V, E_{i\less j})$, where $E_{i\less j}:=\{A\in E_i\less E_j:j\notin A\}$. Let $G-i$ be the $r$-graph obtained from $G$ by deleting vertex $i$ and the edges containing $i$. A hypergraph $G=(V,E)$ is said to {\it cover} a vertex pair $\{i,j\}$ if there exists an edge $e\in E$ with $\{i,j\}\subseteq e$. $G$ is said to {\it cover pairs} if it covers every pair $\{i,j\}\subseteq V^{(2)}$.
\begin{lem}[\cite{FF1989,T2017}]\label{Hmr}
Suppose $G\subseteq [n]^{(r)}$ and $\overset{\rightarrow}{x}=(x_1,\dots,x_n)$ is a legal weighting. For all $1\leq i<j\leq n$, we have
\begin{description}
  \item[(i)]  Suppose that $G$ does not cover the pair $\{i,j\}$. Then $\lambda(G)\leq\max\{\lambda(G-i),\lambda(G-j)\}$. In particular, $\lambda(G)\leq\lambda([n-1]^{(r)})$.
  \item[(ii)] Suppose that $m,t\in \mathbb{N}$ satisfy ${{t-1}\choose r}\leq m\leq {t\choose r}-{{t-2}\choose{r-2}}$. Then
  $$\lambda(C_{r,m})=\lambda([t-1]^{(r)})=\frac{1}{(t-1)^r}{{t-1}\choose r}.$$
  \item[(iii)] $w(G_i,\overset{\rightarrow}{x})\leq (1-x_i)^{r-1} \lambda(G_i)$ for any $i\in [n]$.
\end{description}
\end{lem}

% \begin{lem}[\cite{MS1965}]\label{GiGj}
% Let $G=(V,E)$ be an $r$-graph and $\overset{\rightarrow}{x}=(x_1,\dots,x_n)$
% an optimal legal weighting for $G$ with $k\leq n$ nonzero weights $x_1,\dots,x_k$. Then, for every $\{i,j\}\in[k]^{(2)}$,
% \begin{description}
%   \item[(i)] $w(G_i,\overset{\rightarrow}{x})=w(G_j,\overset{\rightarrow}{x})$;
%   \item[(ii)] there is an edge in $E$ containing both $i$ and $j$.
% \end{description}
% \end{lem}
\begin{defn}[\cite{B1986}]
Let $E\subset \mathbb{N}^{(r)}$, $e\in E$ and $i,j\in \mathbb{N}$ with $i<j$. Then define
\begin{equation*}
L_{ij}(e)=\left\{\begin{array}{cc}
                             (e\less\{j\})\cup \{i\}, & \text{if}~i\notin e ~\text{and}~ j\in e;\\
                             e, & \text{otherwise},
                           \end{array}\right.
\end{equation*}
and $$\mathscr{C}_{ij}(E)=\{L_{ij}(e):e\in E\}\cup \{e:e,L_{ij}(e)\in E\}.$$
We say that $E$ is {\it left-compressed} if $\mathscr{C}_{ij}(E)=E$ for every $1\leq i<j$.
\end{defn}

From now on, suppose that ${{t-1}\choose r}\leq m< {t\choose r}$ for some integer $t$. Let $G$ be a graph with $e(G)=m$ which satisfies $\lambda(G)=\lambda_r(m)$ and let $\overset{\rightarrow}{x}$ be a legal weighting attaining the Lagrangian of $G$. Without loss of generality, we can assume $x_i\geq x_j$ for all $i<j$ and $\overset{\rightarrow}{x}$ has the minimum possible number of non-zero entries, and let $T$ be this number.

Suppose that $G$ achieves a strictly larger Lagrangian than $C_{r,m}$. Then we have
%\begin{equation}\label{suppose}
$$\lambda(G)>\frac{1}{(t-1)^r}{{t-1}\choose r},$$
 %\end{equation}
 which in turn implies $T\geq t$, otherwise $\lambda(G)\leq \lambda([t-1]^r)$.

 \begin{lem}[\cite{FF1989,FR1984,T2002}]\label{four}
 Let $G,T$ and $\overset{\rightarrow}{x}$ be as defined above. Then
\begin{description}
\item[(i)] $G$ can be assumed to be left-compressed and to cover pairs.
\item[(ii)] For all $1\leq i\leq T$ we have  $$w(G_i,\overset{\rightarrow}{x})=r\lambda(G_i).$$
\item[(iii)] For all $1\leq i<j\leq T$ we have $$(x_i-x_j)w(G_{i,j},\overset{\rightarrow}{x})=w(G_{i\less j},\overset{\rightarrow}{x}).$$
%\item[(iv)] If $T=t$, then $x_1<\frac{1}{t-r+1}$  and $x_1<\frac{k+1}{k}x_{t-(k+1)r}$ for $1\leq k\leq \frac t r-1$.
\end{description}
\end{lem}

\begin{lem}[\cite{T2017}] \label{x1}
If $T=t$, then $x_1<\frac{1}{t-r+1}$  and $x_1<\frac{k+1}{k}x_{t-(k+1)r}$ for $1\leq k\leq \frac t r-1$.
\end{lem}

\subsection{The case $r=3$}
Let $r=3$, ${t-1\choose 3}\leq m <{t\choose 3}$ for some integer $t$, and $G$ be a $3$-graph with $m$ edges so that $\lambda(G)=\lambda_3(m)>\lambda([t-1]^3)$. Let $\vec x=(x_1,\ldots, x_n)$ be an optimal legal weighting for $G$ that uses exactly $T$ nonzero weights (i.e., $x_1\geq \cdots \geq x_T> x_{T+1}=\cdots =x_n=0$).
Talbot (\cite{T2002}, Inequality (2.2)) proved that the number of edges in $G$ must satisfy
$$m\geq {T-1\choose 3}+{T-2\choose 2} -(T-2).$$
Since $m<{t\choose 3}$ and $T\geq t$, it implies $T=t$. Thus,
\begin{lem}[{\cite{T2002}}] \label{T=t}
  $G$ must have support on exactly $t$ vertices, i.e., $T=t$.
\end{lem}

% \begin{lem}\label{jkl}
% For any $j,k,l\in [t]$, we have
% \begin{align}\label{diff}
% x_1x_{t-1}x_t-x_jx_kx_l\leq (x_1-x_j-x_k-x_l+x_{t-1}+x_t)x_{t-1}x_t.
% \end{align}
% \end{lem}
% \pf To prove \eqref{diff}, it suffices to verify the following inequality,
% \begin{align}\label{difference}
% %&(x_1-x_j-x_k-x_l+x_{t-1}+x_t)x_{t-1}x_t-(x_1x_{t-1}x_t-x_jx_kx_l)\notag\\
% x_j(x_kx_l-x_{t-1}x_t)-(x_k+x_l-x_{t-1}-x_t)x_{t-1}x_t\geq0.
% \end{align}
% However
% \begin{align*}
% x_j(x_kx_l-x_{t-1}x_t)&=x_j(x_kx_l-x_kx_t+x_kx_t-x_{t-1}x_t)\\
% %&=x_j(x_k(x_l-x_t)+x_t(x_k-x_{t-1}))\\
% &=(x_l-x_t)x_jx_k+(x_k-x_{t-1})x_jx_t\\
% &\geq(x_k+x_l-x_{t-1}-x_t)x_{t-1}x_t.
% \end{align*}
% Hence Inequality \eqref{difference} holds.
% \proofsquare

Lemmas \ref{x1} and \ref{T=t} imply the following inequality:
\begin{align}
  x_1 &< \frac{1}{t-2}. \label{eq:x1}
  %x_1 &<\frac{k+1}{k}x_{t-3(k+1)}. \label{eq:xk}
\end{align}
We have the following lemmas:
\begin{lem}
  For any $k\in[t-1]$, we have
  \begin{equation}\label{eq:xk}
   x_{t-k}>\frac{k-1}{k+1}x_1.
  \end{equation}
\end{lem}
\pf Observe that \begin{align*}
1&=x_1+\dots+x_{t-k-1}+x_{t-k}+\dots+x_t\\
&< (t-k-1)x_1+(k+1)x_{t-k}\\
&< \frac{t-k-1}{t-2}+(k+1)x_{t-k}.~~~~~~(by~\eqref{eq:x1})
\end{align*}
Solving $x_{t-k}$,  we get
$$x_{t-k}>\frac{k-1}{k+1}\cdot\frac{1}{t-2}>\frac{k-1}{k+1}x_1.$$
\proofsquare

\begin{lem}
  For any subset $S\subseteq [t]$, we have
  \begin{equation}
    \label{eq:sumxk}
    \sum_{i\in S}(x_1-x_i)<2x_1.
  \end{equation}
\end{lem}

\pf It is trivial when $|S|\leq 2$. We can assume $|S|>2$. We will prove it by contradiction.
Suppose that there is $S=\{i_1,i_2,\ldots, i_k\}$ of $k$ distinct elements
such that $x_{i_1}+x_{i_2}+\dots+x_{i_k}\leq(k-2)x_1$.
Then
\begin{align*}
1&=x_1+x_2+\dots+x_t\\
&\leq x_{i_1}+x_{i_2}+\dots+x_{i_k}+(t-k)x_1\\
&\leq (k-2)x_1+(t-k)x_1\\
 &=(t-2)x_1\\
&<1.
\end{align*}
Contradiction.
%Then $x_1\geq\frac{1}{t-2}$,
%contrary to Lemma \ref{four} (iv).
\proofsquare

\section{Proof of Theorem \ref{mainthm}}

{\bf Proof of Theorem \ref{mainthm}:} Write $m={t\choose 3}-l$ where $0< l\leq {t-1\choose 2}$.
Let $m'=m+\eta$, where $\eta:=\lceil 4t^{2/3}\rceil$.

We claim for all $t\geq 8$
\begin{equation}
  \label{eq:m}
  \lambda_3(m)\leq \lambda(C_{m'}).
\end{equation}

Without loss of generality,
we can assume $\lambda_3(m)>\lambda([t-1]^{(3)})$. Otherwise, we have
$$\lambda_3(m)\leq \lambda([t-1]^{(3)})\leq \lambda(C_{m'}).$$

When $l\leq \eta $, then $m'\geq {t\choose 3}$. We have
$$\lambda_3(m)\leq \lambda_3\left({t\choose 3}\right)=\frac{{t\choose 3}}{t^3}\leq \lambda(C_{m'}).$$

We can assume $l>\eta$. Let $l'=l-\eta\geq 1$.

Let $G=(V,E)$ be a $3$-graph with $m$ edges satisfying $\lambda(G)=\lambda_3(m)$ and
Let $\vec x=(x_1,\ldots, x_n)$ be an optimal legal weighting for $G$ that uses exactly $T$ nonzero weights (i.e., $x_1\geq \cdots \geq x_T> x_{T+1}=\cdots =x_n=0$). By Lemma \ref{T=t}, we have $T=t$. By considering the induced sub-hypergraph on first $t$ vertices, we can assume $G$ has exactly $t$ vertices, at most $m$ edges,
%{\color{red} and Lemma \ref{four}(i), we can assume $G$ has exactly $t$ vertices and $m$ edges},
and $\lambda(G)=\lambda_3(m)$. In addition, we may assume $G$ is left-compressed
by Lemma \ref{four}(i).
Define $b=\max\{i: \{i,t-1,t\} \in E\}$, we have $E_i=([t]\less i)^{(2)}$ for $i\in [b]$. So $x_1=x_2=\dots=x_b$ by Lemma \ref{four} (iii).

By the definition of $b$, we have
$$E^c_{t-1,t}=\{i \colon b+1\leq i\leq t-2\}.$$
In particular,
$G$ is a subgraph of $C_{{t\choose 3}- (t-2-b)}$. If $t-2-b\geq l'$, then $G$ is also a subgraph of $C_{m'}$. This implies
$$\lambda_3(m)=\lambda(G)\leq \lambda(C_{m'}),$$
and we are done. So we may assume $t-2-b<l'$. Let $l''= b+ \min\{l'-(t-2), 0\}$.

Let $B\subseteq E^c\setminus \{\{i,t-1,t\} \colon i\in E^c_{t-1,t}\}$ be any set of $l''+\eta$ non-edges.
This is possible since $G$ has at least $l$ non-edges and
$$l=l'+\eta\geq l''+\eta+(t-2-b).$$
Let $G'$ be a $3$-graph obtained from $G$ by deleting all edges in $\{\{b+1-i, t-1,t\}\colon 1\leq i \leq l''\}$
and adding all triples in $B$ as edges. Then $G'$ has at most $m+\eta= m'$ edges.
The main proof is to show the following inequality:
\begin{equation}
  \label{eq:main}
  w(G,\vec x)\leq w(G',\vec x).
\end{equation}

Let $s=\max\{i: \{t-i-1,t-i \}\in E^c_t\}$ and $S=\{t-s,t-s+1,\dots,t-1,t\}$. By the choice of $s$, we know $\{t-s-1,t-s,t\}\in E^c$
but $\{t-s-2,t-s-1,t\}\in E$.
%So $\{b+1,b+2,\dots,t-s-1\}\times \{t-s+1,t-s+2,\dots,t\}\subseteq E_{t-s}$. \bigskip

\noindent {\bf Claim\refstepcounter{counter}\label{eS}  \arabic{counter}.} For any $e\in E^c$, we have $|e\cap S|\geq 2$.

\pf Suppose $e=\{i,j,k\}\in E^c$ with $i<j<k$ and $|e\cap S|\leq 1$. We must have $i,j\notin S$. Since $E$ is left-compressed, $\{i,j,t\}\in E^c$. Then $\{j-1,j,t\}\in E^c$, contrary to the choice of $s$.
\proofsquare

%By Lemma \ref{eS}, we have $E^c\subseteq \{ijt|ij\in E^c_t\}\bigcup(\cup_{\{i,j\}\subseteq S^{(2)}}F_{ij})$.

% \noindent {\bf Claim\refstepcounter{counter}\label{1t-s}  \arabic{counter}.}
%   $x_1<\frac{s}{s-3}x_{t-s}$  and  $x_1<\frac{s-1}{s-4}x_{t-s+1}$

% \pf Since $x_1<\frac{k+1}{k}x_{t-3(k+1)}$ for $1\leq k\leq \frac t r-1$ by Lemma \ref{four} (iv), we can get the two bounds by making $k=\frac{s-3}3$ and $k=\frac{s-4}3$, respectively.
% \proofsquare

We may assume
 \begin{align}\label{eq:x1xt-1xt}
x_1x_{t-1}x_t-x_{t-s-1}x_{t-s}x_t\geq 0.
\end{align}
Otherwise by replacing the edge $\{1,t-1,t\}$ with the non-edge $\{t-s-1, t-s, t\}$,
we get another 3-graph with the same number of edges whose Lagrangian is strictly greater than the Lagrangian of $G$.
%By Lemma \ref{jkl}, we have $$x_1x_{t-1}x_t-x_{t-s-1}x_{t-s}x_t\leq (x_1-x_{t-s-1}-x_{t-s}+x_{t-1})x_{t-1}x_t.$$
% Thus
% \begin{align}\label{eq:x1}
%x_1-x_{t-s-1}-x_{t-s}+x_{t-1}\geq 0.
%\end{align}

 Combining Inequalities \eqref{eq:x1xt-1xt} and \eqref{eq:xk}, we get
 \begin{align}\label{eq:x1xt-1}
   x_{t-1}\geq\frac{x_{t-s-1}x_{t-s}}{x_1}> \frac{s(s-1)}{(s+2)(s+1)}x_{1}.
   \end{align}

For any $\{j,k\}\subseteq S^{(2)}$ with $j<k$, let $F_{jk}=\{\{i,j,k\}|i\in E^c_{jk}~\text{and}~ i<j\}$.
By Claim \ref{eS}, we have $$E^c= \bigcup\limits_{\{i,j\}\subseteq S^{(2)}}F_{ij}.$$

Now we will prove Inequality \eqref{eq:main}. We divide it into two cases.

\noindent
{\bf Case 1:} $s^2+s< \eta$. We have

\begin{align*}
  w(G',\vec x) &= w(G, \vec x) - l''x_1x_{t-1}x_{t}+ \sum_{\{i,j,k\}\in B} x_ix_jx_k\\
                &=w(G, \vec x)+\eta x_1x_{t-1}x_t -\sum_{\{i,j,k\}\in B}(x_1x_{t-1}x_t-x_ix_jx_k)\\
                &>w(G, \vec x)+\eta x_1x_{t-1}x_t - \sum_{\{i,j,k\}\in B} (x_1-x_i)x_{t-1}x_t\\
                 &\geq w(G, \vec x)+\eta x_1x_{t-1}x_t - x_{t-1}x_t\sum_{\{j,k\}\in S^{(2)}} \sum_{i\colon \{i,j,k\}\in F_{jk}} (x_1-x_i)\\
                &\geq w(G, \vec x)+\eta x_1x_{t-1}x_t - x_{t-1}x_t\sum_{\{j,k\}\in S^{(2)}} 2x_1   ~~~~~~~~~~~~~~~~~~~(by~ \eqref{eq:sumxk})\\
                &=w(G, \vec x)+ (\eta-s^2-s)x_1x_{t-1}x_t\\
                &>w(G, \vec x).
\end{align*}

\noindent
{\bf Case 2:} $s^2+s\geq \eta$. We have
\begin{equation}
  \label{eq:s}
s\geq \frac{\sqrt{1+4\eta}-1}{2}>\sqrt{\eta}-\frac{1}{2}.
\end{equation}
Since $\eta=\lceil 4t^{2/3}\rceil$ and $t\geq 8$, by Inequality \eqref{eq:s}, we have
\begin{equation}
  \label{eq:seta}
s(s-1)\eta>(4s+2)(t-2)\geq (4s+2)l''.
\end{equation}
\begin{align*}
  w(G',\vec x) &= w(G, \vec x) - l''x_1x_{t-1}x_{t}+ \sum_{\{i,j,k\}\in B} x_ix_jx_k\\
                &>w(G, \vec x)-l'' x_1x_{t-1}x_t + (l''+\eta) x_{t-1}^2x_t\\
                &=w(G, \vec x)+ x_{t-1}x_t\left((l''+\eta) x_{t-1}-l'' x_1\right)\\
                &\geq w(G, \vec x)+ x_{t-1}x_t\left((l''+\eta) \frac{s(s-1)}{(s+2)(s+1)} x_1- l'' x_1\right)~~~~~~(by~\eqref{eq:x1xt-1})\\
                &=w(G, \vec x)+ \frac{1}{(s+2)(s+1)}x_1x_{t-1}x_t\left(s(s-1)\eta -(4s+2)l''\right)\\
               &> w(G, \vec x). ~~~~~~~~~~~~~~~~~~~~~~~~~~~~~~~~~~~~~~~~~~~~~
                 ~~~~~~~~~~~~~~~~~~~(by~\eqref{eq:seta})
\end{align*}
Therefore, Inequality \eqref{eq:main} holds in any circumstances. If $l' \leq t-2$, then $G'$ is a subgraph of $C_{m'}$, else $G'$ is a subgraph of
$C_{{t\choose3}-(t-2)}$.
Inequality \eqref{eq:m} follows from Inequality \eqref{eq:main} by a sequence of inequalities:\\

if $l'\leq t-2$, then
$$\lambda_3(m)=w(G,\vec x)\leq w(G',\vec x)\leq \lambda(G')\leq \lambda(C_{m'}),$$
else
$$\lambda_3(m)=w(G,\vec x)\leq w(G',\vec x)\leq \lambda(G')\leq\lambda\left(C_{{t\choose3}-(t-2)}\right)=\lambda([t-1]^{(3)})  \leq \lambda(C_{m'}).$$
Finally we can choose a constant $c$ large enough so that the following two conditions hold:
\begin{itemize}
\item $cm^{2/9}>4\lceil t^{2/3} \rceil$ for all $t\geq8$,
\item and $cm^{2/9}>{t-1 \choose 2}$ for $1\leq t\leq 8.$
\end{itemize}
When $t\geq 8$, we have
$$\lambda_3(m)\leq \lambda(C_{m'}) \leq \lambda(C_{m+cm^{2/9}}).$$
When $1\leq t\leq 8$, we have
$$m+cm^{2/9}>{t-1\choose 3}+{t-1\choose 2}={t\choose 3}.$$
We have
$$\lambda_3(m)\leq \lambda_3\left({t\choose 3}\right)=\frac{{t\choose 3}}{t^3}=
\lambda\left(C_{{t\choose 3}}\right)\leq \lambda(C_{m+cm^{2/9}}).$$

This completes the proof of Theorem~\ref{mainthm}.
\proofsquare

{\bf Remark:} Actually Inquality \eqref{eq:seta} only requires $\eta=c(l'')^{2/3}$.
When $m$ is closed to ${t\choose 3}$, we can get a better bound.
Let $m={t\choose 3}-l$ where $0<l<(t-2)+ct^{2/3}$.
Then we have
$$\lambda_3(m)\leq \lambda(C_{m+cl^{2/3}}).$$

\end{document}